\documentclass[12pt, reqno]{amsart}
\usepackage{amsmath, amsthm, amscd, amsfonts, amssymb, graphicx, color}
\usepackage[bookmarksnumbered, colorlinks, plainpages]{hyperref}

\theoremstyle{definition}

\theoremstyle{remark}

\numberwithin{equation}{section}

\input{mathrsfs.sty}
\begin{document}
\setcounter{page}{1}

\noindent\parbox{2.95cm}{\includegraphics*[keepaspectratio=true,scale=0.125]{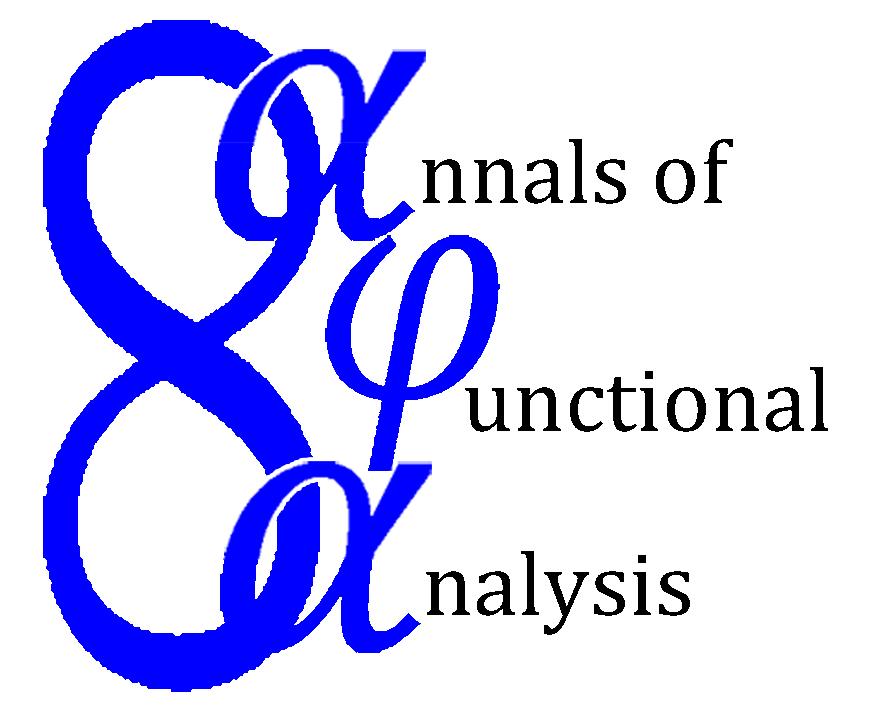}}
\noindent\parbox{4.85in}{\hspace{0.1mm}\\[1.5cm]\noindent
\qquad Ann. Funct. Anal. x (xxxx), no. x, xx--xx\\
{\footnotesize \qquad \textsc{\textbf{$\mathscr{A}$}nnals of
\textbf{$\mathscr{F}$}unctional
\textbf{$\mathscr{A}$}nalysis}\\
\qquad ISSN: 2008-8752 (electronic)\\
\qquad URL:\textcolor[rgb]{0.00,0.00,0.99}{www.emis.de/journals/AFA/}}\\[.5in]}

\title{An interview with Tsuyoshi Ando}

\author[R.A. Brualdi, M.S. Moslehian]{Richard A. Brualdi$^1$ and Mohammad Sal Moslehian$^*$$^2$}

\address{$^1$ Math Dept., UW-Madison, 480 Lincoln Drive, Madison, WI 53706, USA.}
\email{\textcolor[rgb]{0.00,0.00,0.84}{mailto:brualdi@math.wisc.edu}}
 \address{$^2$ Tusi Mathematical Research Group (TMRG), P.O. Box 1113, Mashhad 91775, Iran}
\email{\textcolor[rgb]{0.00,0.00,0.84}{moslehian@member.ams.org}}

\begin{abstract}
In celebration of the distinguished achievements of Professor Tsuyoshi Ando in matrix analysis and operator theory, we conducted an interview with him via email. This paper presents Professor Ando's responses to several questions we gave him regarding his education and life as a mathematician.
\end{abstract}



\date{Received: 15 October 2013; Accepted: 30 January 2014 .
\newline \indent $^{*}$ Corresponding author}
\maketitle

\section{Introduction}

Professor Tsuyoshi Ando has made fundamental contributions in many different parts of matrix analysis and operator theory and has written extensively on these subjects. Several theorems, concepts, theories and inequalities carry his name: ``Krein--Ando theorem", ``Ando--Hiai inequality", ``Ando--Li--Mathias
geometric mean", ``Ando--Douglas type theorem in Riesz spaces",
``Hua--Marcus--Bellman--Ando inequalities on contractive matrices",
``Kubo--Ando theory of operator means", ``Ando--Krieger theorem", ``The
central Ando dilation ", ``Amemiya--Ando theorem", ``Ando's inequality for
positive linear maps". These are evidence of his substantial influence on the development of these subjects.

He has served as a member of the editorial board of several prestigious journals such as {\it Linear Algebra and its Applications}, of which he is now a Distinguished Editor, {\it Operators and Matrices}, and {\it Positivity}. He has been given several awards:

\begin{itemize}
\item The 2002 Hans Schneider Prize in Linear Algebra (The International Linear Algebra Society);
\item The 2005 B\'ela Sz\"okefalvi Nagy Prize (Acta Scientiarum Mathematicarum, Hungary);
\item The 2011 Second Order of the Scared Treasure (Japan Government).
\end{itemize}

In addition, several issues of journals have been dedicated to him as follows:

\begin{itemize}
\item {\it Operator Theory, Advances and Applications} vol. $62$, Contributions to Operator Theory and its Applications, 1993.
\item {\it Linear Algebra and its Applications}, volume 341 (2002).
\end{itemize}

Although the interviewers have visited Professor Ando several times and have been at many conferences with him, the following interview is based on some communications via email.
We are grateful to Professor Ando for granting this interview and for his substantial informative responses.
\section{Interview}


$\bullet$ \textbf{Tell us about your childhood and your education, beginning with
elementary school. Tell us about your mathematics education in
college including graduate work.}\\

$\star$
In 1932 I was born in Sapporo, the local capital of Hokkaido, a northern
island of Japan, and have grown up there and will die probably there.

Sapporo is a snowy city. When I was born, the population of Sapporo was about 200,000, and now it is about 1.8 million.

I entered an elementary school in 1938. Elementary school was compulsory education of 6 years. Usual school children were classified
into boys' classes and girls' classes separately. Only children of poor
health were gathered into common co-education classes.\\

\centerline{\includegraphics[keepaspectratio=true,scale=0.125]{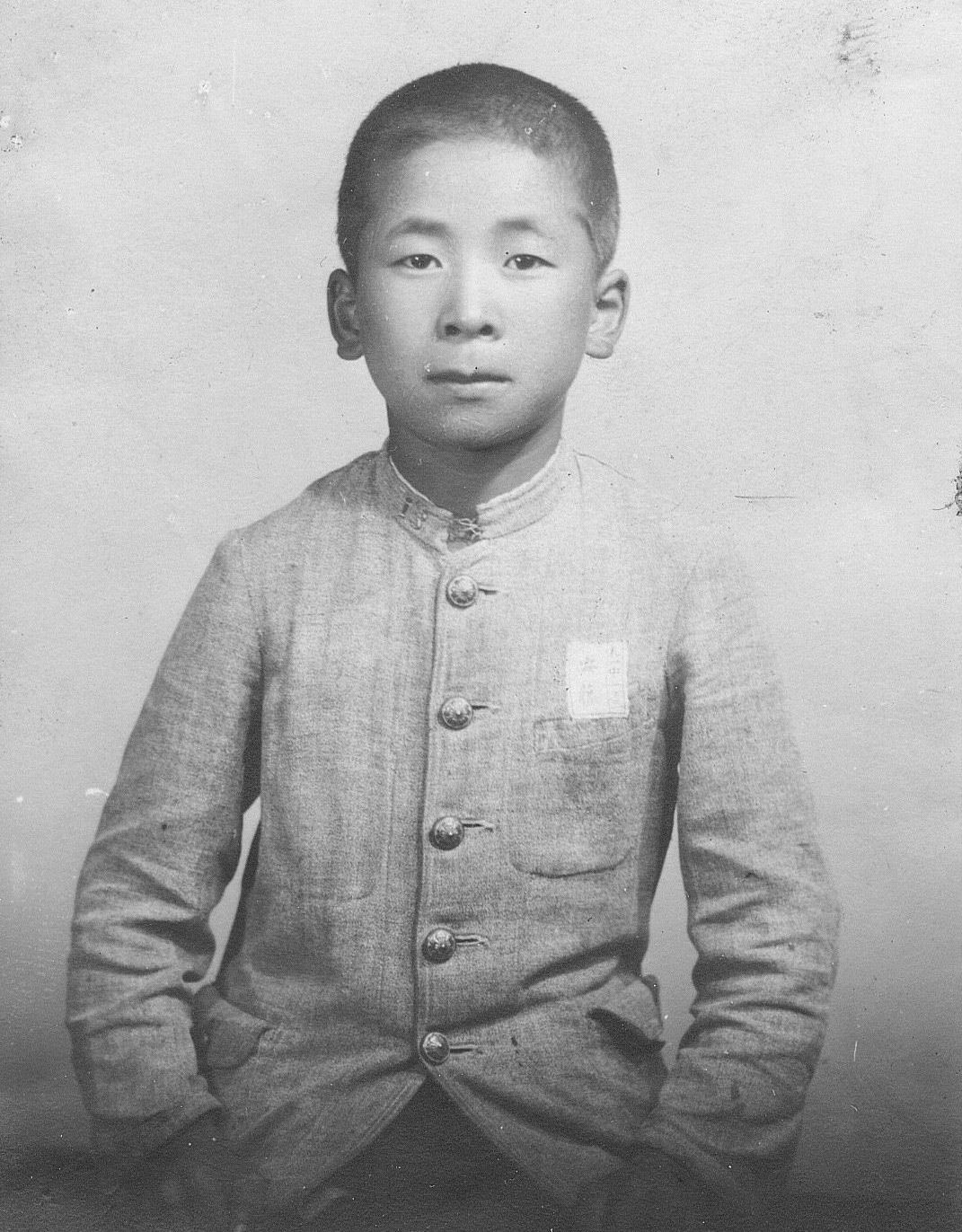}}
\vspace{.05in} \centerline{(T. Ando at the 6-th grade of national school, 1943)}
\vspace{0.5cm}

At that time Japan was already in a war situation with China. At the
end of 1941 Japan opened a war against USA, England and other countries.
Every elementary school changed its name to ``national school".

Until the middle of 1942 Japan seemed to maintain superior position over its
enemies, but thereafter it was always in inferior position.

In April of 1944 I entered a middle school. No special education for
war was given. But boys of the upper classes were sent to farm villages
to help with the harvest and/or cultivation for months.

Someone has said that during the war time English education was suspended in middle school education, but this is not true. We listened to English classes
even though class hours were not many.

In April of 1945, I entered a military cadet school in Sendai, a preparatory school for Military Academy. I really wanted
to become a military officer. Except for some outside camping, special military lessons were not given to the first year boys.

Though heavy bombing to many Japanese cities by U.S. Air Force had started since 1944,
Sendai suffered a serious bombing in July 1945 for the first time. The city of Sendai burned out in one night. The school began to prepare for evacuation to some village in
Fukushima. Just on the day when the first group of evacuees started, Japan surrendered. The cadet school was closed. I returned back to Sapporo and back to the same middle school.

Education in middle school after the war was in confusion. Much was not taught, even in mathematics. I liked mathematics mostly because the answers in examinations of mathematics were always unique and I could get good notes.

Until the end of the war, the education system in Japan consisted of a
6-year elementary school, a 5-year middle school, a 3-year national high school,
and a 3-year university. Here national high school is different from the
present day high school. Its education level corresponded to that of the
first two years of general education curricula of today's 4-year university.

After the war, following the suggestion of the
USA, Japan changed its education system as follows: a 6-year elementary
school, a 3-year middle school, a 3-year high school, and a 4-year university.
We were in a transitional period.

In the old system of education one could apply for high school at the end of the
4-th year of middle school. Finishing the 4-th year of middle school,
I entered the preparatory school of
Hokkaido University (which was equivalent to a high school).

Hokkaido University was the only university of the old system in Hokkaido.
But the following year the new system began and we were shifted to the Hokkaido University of new system as first year students.

In the new system basic curricula for freshmen seemed not to be well organized at first. Each lecturer taught some things according his interest. For instance, I listened to courses on elementary dynamics from three different lecturers, but nothing about electromagnetism.

After one and a half years of general education each student had to choose his major subject. My brother, who was in the same grade as mine, chose electrical engineering.
I though it was time for me to go in a different direction and went to science.
I thought of either physics or mathematics. And finally I decided to go to mathematics.

I could not draw my future plan. Most of graduates in mathematics had found their ways
as school teachers. As the economy of Japan was not active yet, there was no possibility
for graduates of mathematics to find jobs in industry. To get an academic position seemed
almost hopeless.\\

\centerline{\includegraphics[keepaspectratio=true,scale=0.125]{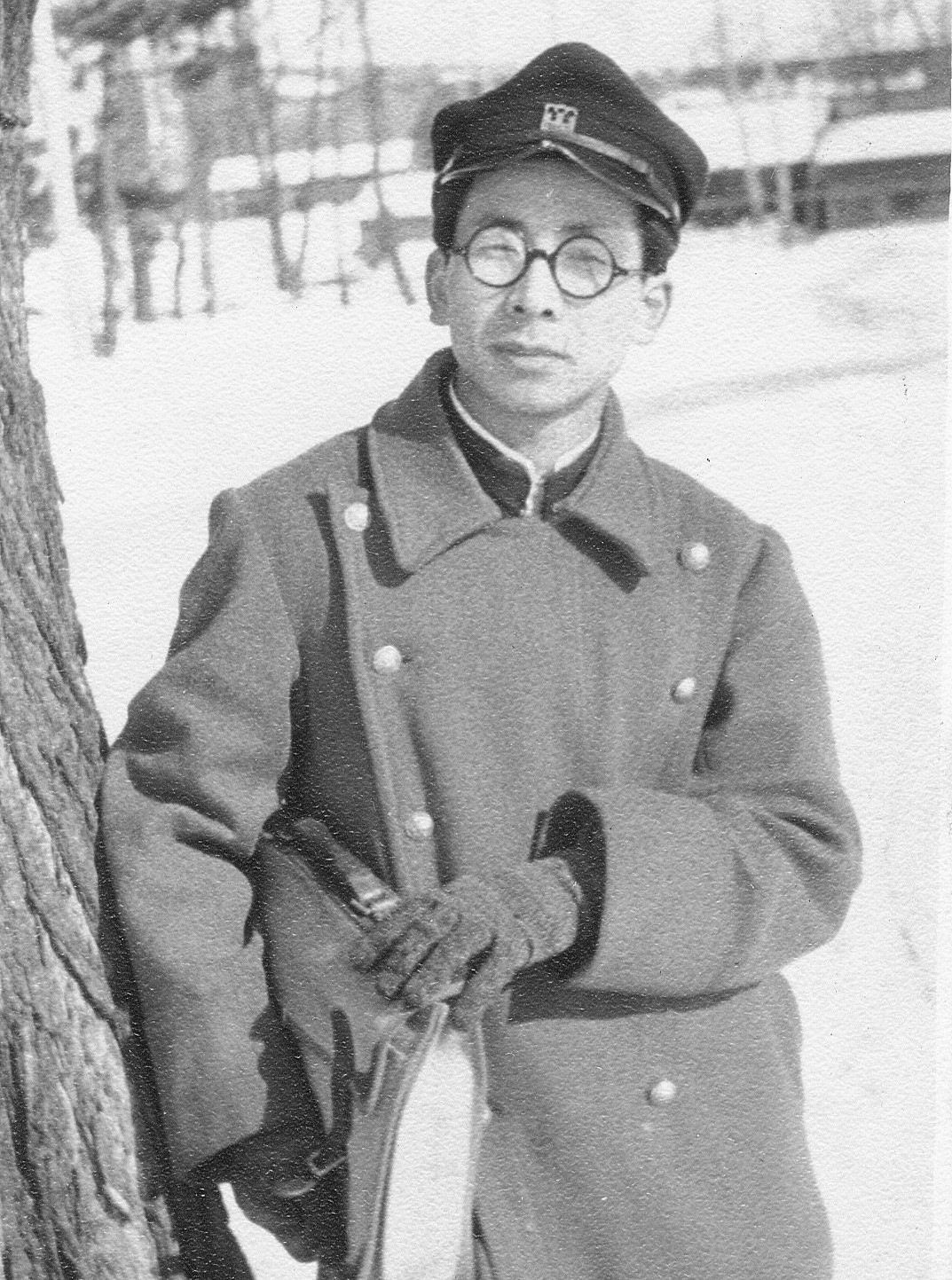}}
\vspace{.05in} \centerline{(T. Ando in the university year 6, 1952)}
\vspace{0.5cm}

Education in the Department of Mathematics was also not satisfactory for different reasons. Originally the Department had four full professors (in the German sense).
Those professors were not on good terms, and three of them moved to other universities
and so did even most assistant professors. The atmosphere in the Mathematics Department was not stimulating. Even basic courses were delivered by
part time visiting lecturers. Then Dr.~Hidegoro Nakano was invited
as professor of functional analysis from Tokyo University.\\

\centerline{\includegraphics[keepaspectratio=true,scale=0.125]{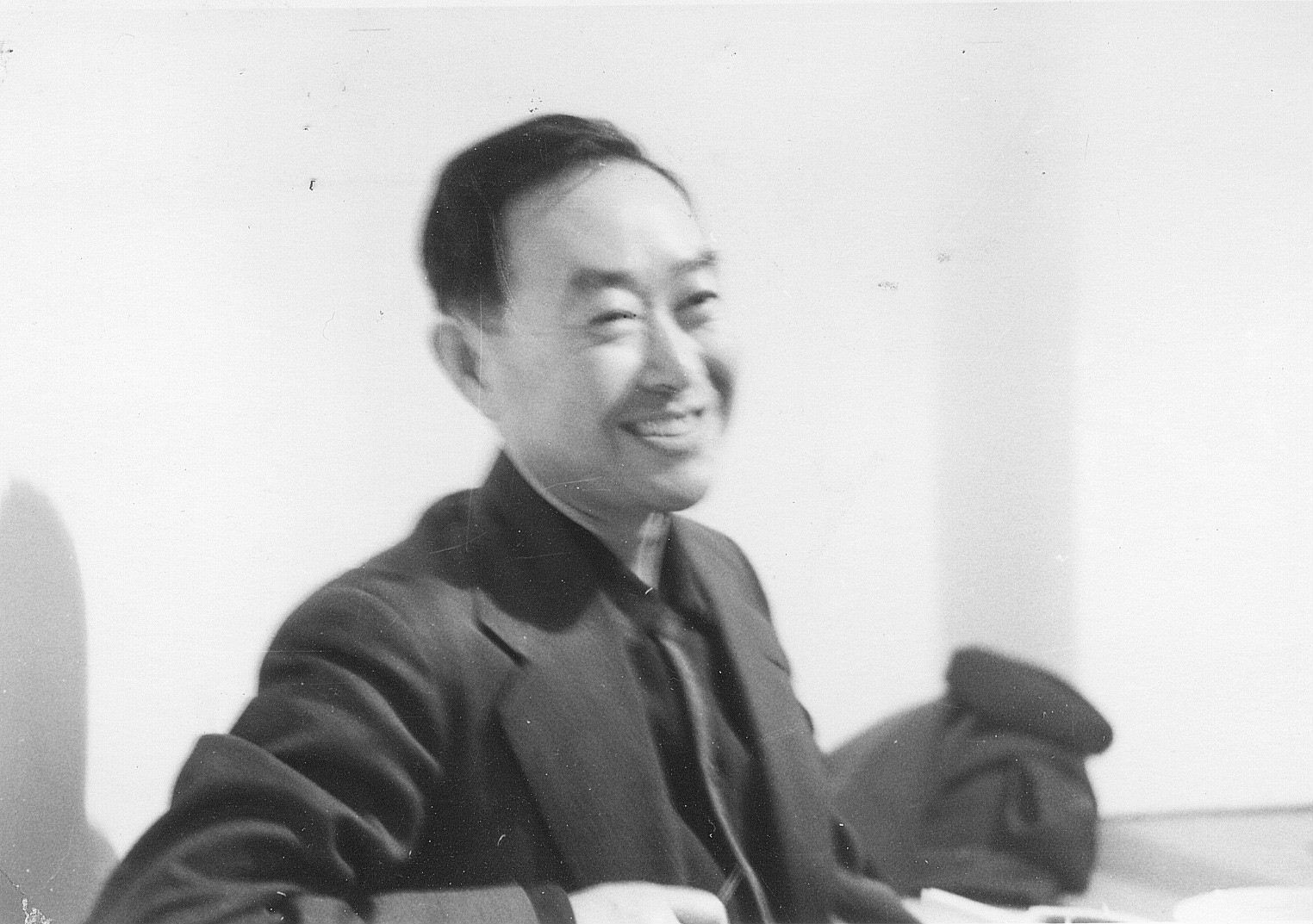}}
\vspace{.05in} \centerline{(H. Nakano, around 1955)}
\vspace{0.5cm}

In graduate study, as the (only remaining) professor of geometry was of the old style, it was easy for me to choose
Prof. Nakano as a supervisor. Four students came under his supervision.

Each student had to select a monograph for the elementary seminar (of training).
T. Shimogaki and I chose S. Banach: ``Th\'eorie \' Operations lin\'earires" (1932) while
other two, J. Ishii and M. Sasaki, chose B. Sz.-Nagy: ``Spektraldarstellung linearer Transformationen des Hilbertschen Raumes" (1942).
Banach's book, presenting his original ideas, was not written in expository style and was very difficult to read for a beginner. We could not come to the end of the book mostly because of a long hospital stay of the supervisor.\\

\centerline{\includegraphics[keepaspectratio=true,scale=0.175]{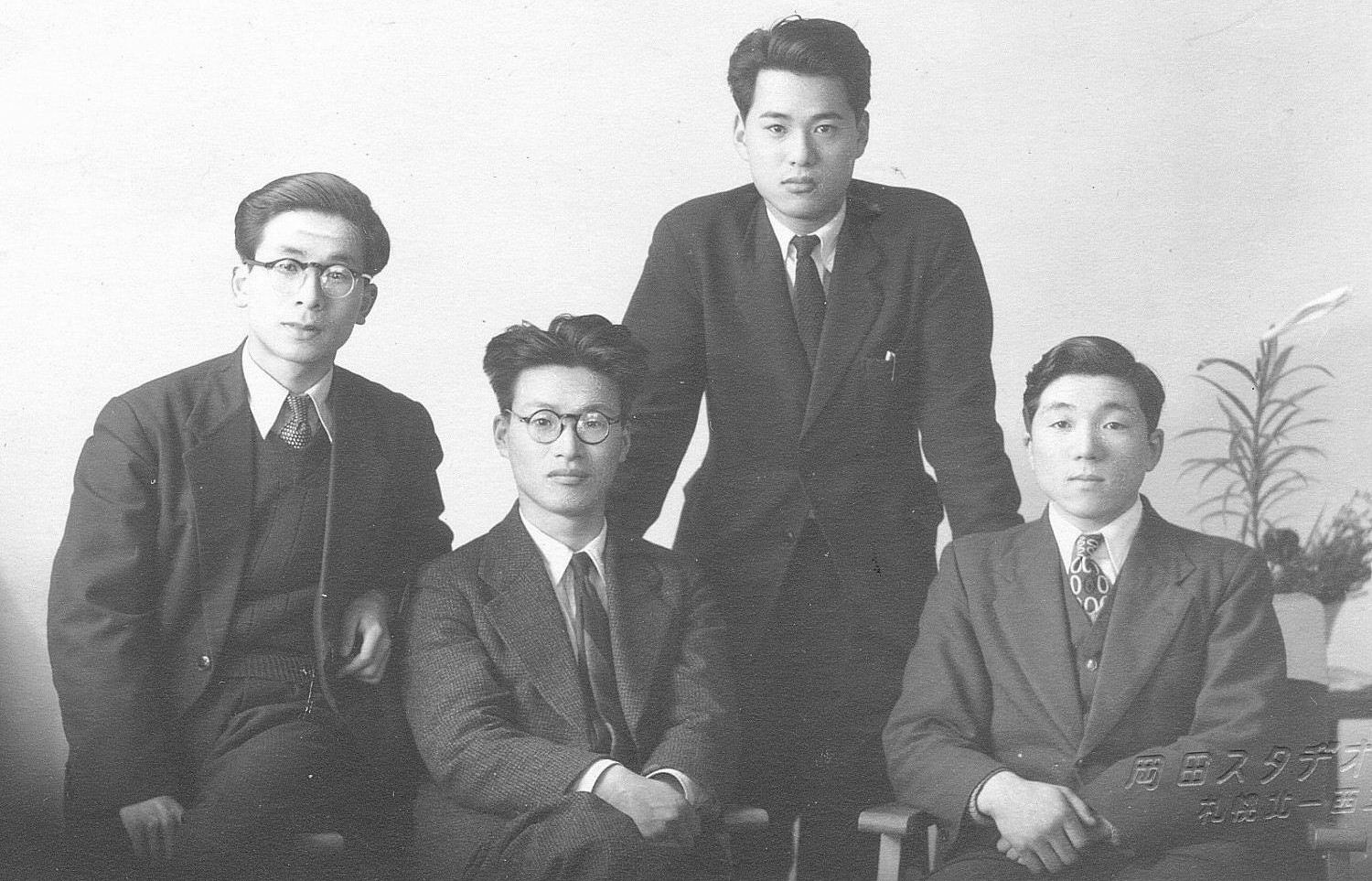}}
\vspace{.05in} \centerline{(T. Shimogaki (back)
    T. Ando (left), J. Ishii (middle), M. Sasaki(right), 1955)}
\vspace{0.5cm}

\noindent Prof. Nakano had already published several monographs:\\
- ``Modulared semi-ordered linear spaces'', 288 pages, 1950;\\
- ``Modern Spectral Theory", 323 pages, 1950;\\
- ``Topology and Linear Topological Spaces" 281 pages, 1951;\\
- ``Spectral Theory in the Hilbert space", 300 pages, 1953.

I read in detail ``Modern Spectral Theory" and was much impressed
by the beauty of his representation theory of order-complete
vector lattices. In ``Modulared Semi-ordered Linear Spaces''
he developed his theory of modulared semi-ordered linear space, the so-called Orlicz--Nakano spaces,
in full details. Everything was written in abstract form
and no concrete examples were presented.

Beside the elementary seminar for training, we joined the
seminar of the research group conducted by Prof.~Nakano. In each session
someone introduced some new results from outside or his own
contribution. I suppose there was no central theme.

When I was asked to make a presentation in the seminar for the first time, I chose a survey of linear or non-linear maps in $C(X)$, for which I had read
very many papers. After my presentation, Prof.~Nakano mentioned
``What is necessary for mathematical research is not wide
knowledge but originality!" Though his comment was true,
I was a little disappointed and thought in secret ``Knowledge is also necessary."

In the custom of Japan at that time, a graduate student had to find the subject
 for his thesis by himself, no hint, no suggestion or no
 help came from the supervisor. We were in desperate atmosphere !

Prof.~Nakano once told us graduate students: ``A graduate student is happy while he has something to read. But once he tries to show his originality then he finds true difficulty of research."

Though the representation theory of Nakano for order-complete vector lattices was quite beautiful, its success was mostly from ``commutativity'' of all order-projections. From this point the space
 of order-bounded linear operators of a vector-lattices
 could not be well treated. Nakano developed his theory of the tensor-product of vector-lattices
 in connection with linear operators on vector-lattices.
 The results was not fruitful because of lack of ``commutativity" in some sense.
 A Perron--Frobenius type theorem
 for positive linear operators was not in the scope of his theory.

 I tried to find a generalization of non-quasi-nilpotency of a (entry-wise) positive
 matrix to a positive linear operator on a Banach lattice.
 I got a generalization under rather strong assumption and wrote a paper \cite{T1},
which became my PhD thesis.

There had been a monumental contribution in this direction by Soviet mathematicians, M.G. Krein and M.A. Rutman:
``Linear operators leaving invariant a cone in a Banach space"
{\it Usphki Matem. Nauk (N.S.)} 3(1948), no. 1 (23), 3--95.
But I could not access the details of the Krein--Rutman paper at that time.

Though my result is sometimes called the Ando--Krieger theorem, the content of the paper was rather modest. I published more than 100 papers since. Many of them are humble, and I imagine they could disappear from the printed journals!

To get an academic position seemed very difficult. As I was the first PhD in the new system and PhD's of the old system were quite few, Prof. Nakano squeezed me into a
position of an academic assistantship of the Mathematics Section of the Research Institute of Applied Electricity, which belonged to the University. The Institute was
established during the war and consisted
of 9 sections. Beside one in mathematics, all the other sections are in experimental sciences;
physics, chemistry, electricity and physiology.

I wondered how long I could stay in the Institute. However, I have stayed in the Institute
until retirement, about 40 years ! I was always isolated, but it was not bad because I could continue to do what I liked.\\

\centerline{\includegraphics[keepaspectratio=true,scale=0.185]{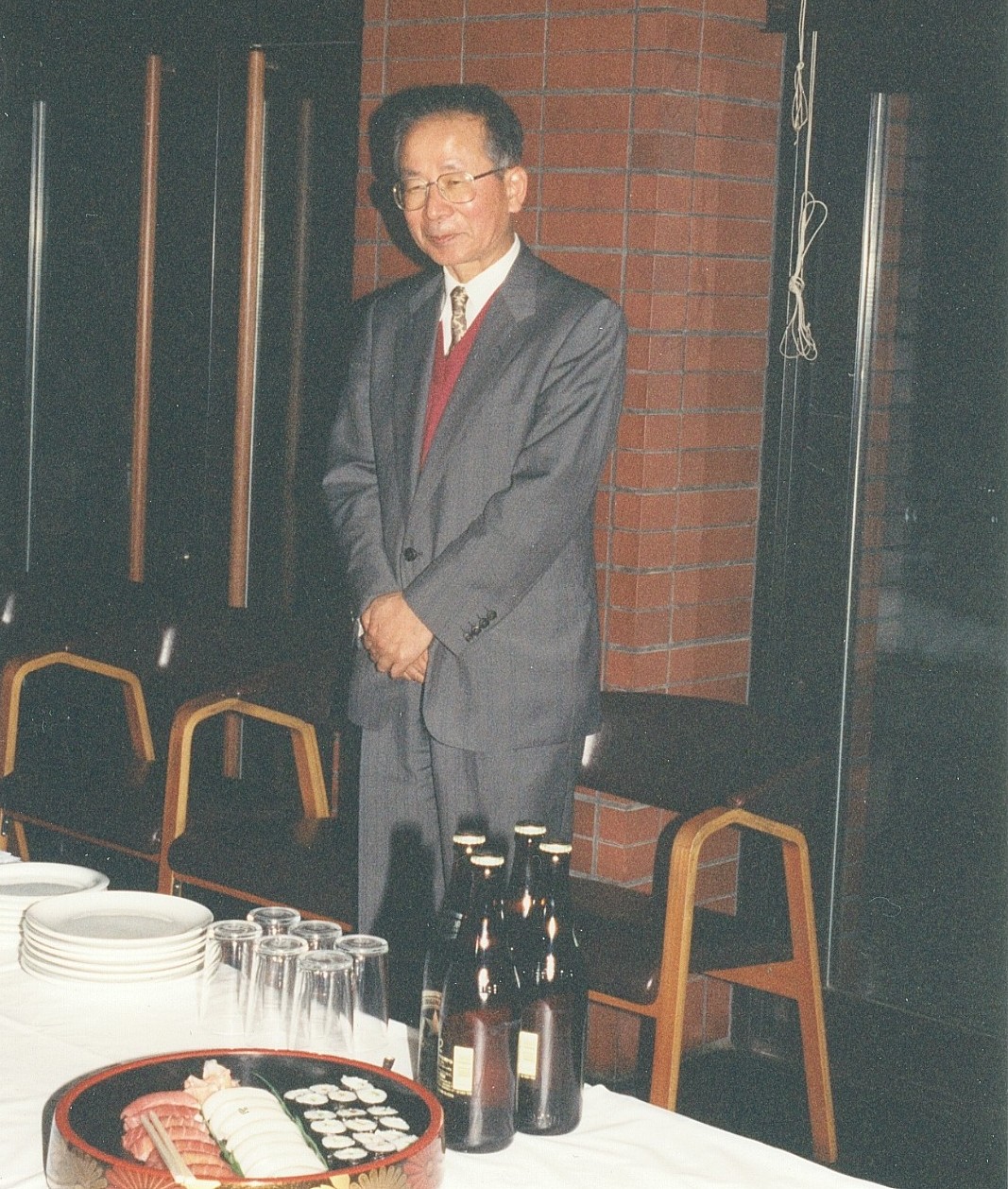}}
\vspace{.05in} \centerline{(T. Ando at a party for retirement, 1995)}

\vspace{.2in} $\bullet$ \textbf{How did your family influence your education? Were either of your
parents educators or professionals? What did they do?}\\

$\star$
My father was a city officer of middle class. I had two brothers. The eldest brother became also a city officer afterwards. There was no academic atmosphere around me.

The other brother was near to me in age. During elementary school
he suffered tuberculosis and had to stay one more year in the same grade.
As a consequence he and I came to be in the same grade. Such situations
were not favorable for me and also for my brother. This brother became
an engineer in industry afterwards.

\vspace{.2in}
$\bullet$ \textbf{When and how did you become interested in mathematics? When did you
realize that you were good in mathematics? Did you ever consider being
something other than a mathematician?}\\

$\star$
In middle school I thought that
my ability was fit more to science than to humanities. As, under the poor situations after the war, experiments by pupils themselves were almost impossible, my preference was
mathematics and physics. I did not think that I was specially talented in mathematics.
In middle school I could not find any teacher inspiring me to continue in mathematics.

I suppose that I could have become an engineer in a company. But as I had an unsociable character, I might be not cooperative with other colleagues and would not be
successful.

Let me tell a short story of my encounter with a small topic of mathematics.
When I was a undergraduate student of mathematics, I bought a newly published book:
Hiraku Toyama ``Theory of Matrices", 1952.

It was quite new in the sense that until that time all books on matrices had had title
``Theory of matrices and determinants", and the stress had been placed on determinants. In this book, Schur's theorem that the Schur product of two positive semi-definite matrices is again positive semidefinite was proved on the basis that
the Schur product is a principal submatrix of their tensor product and the tensor product is naturally positive semidefinite. This approach sounded fresh.

My interest in the Schur product might began this time. Afterwards I wrote several papers concerning the Schur product, including \cite{T16, T19, T23, T26}.

\vspace{.2in} $\bullet$ \textbf{In your opinion what is the best part of being a mathematician? Is
there any bad part?}\\

$\star$
One of the nice points of being a mathematician is that he can do, still by himself alone, everything he wants without any financial help from outside or collaborators. This is really nice after retirement as all information is available now at home directly via
internet.

But there remains serious danger. In the fields of experimental sciences one can say he will do such and such research
aiming at some definite results.
Even if the expected result could not be obtained, he can report that the research along this direction was not successful. But in mathematics if one says he will try to solve such and such a problem and is not successful after a period, he can not say that this direction was not successful. A more talented mathematician may find a way to the goal under the same
or a little changed situations. Therefore mathematicians seem under heavier pressure than people in other experimental sciences.

In this connection there is dangerous
temptation to try to get straight generalizations of established nice results like making
Xerox copies of the original idea. Though this attempt is not bad for training, young mathematicians should try to overcome such temptation gradually. They have to bear the pressure!

\vspace{.2in} $\bullet$ \textbf{You have traveled extensively. What are some of your favorite places and experiences?}\\

$\star$
I attended at a number of conferences, but always for a short period. Long stays were
three times.

The first stay (1962--64) in USA was somewhat strange. As I belonged
not to a department of mathematics but
to a research institute of electronics, where I was almost an only mathematician, the director, a chemist, told me to take a postdoctoral position in a chemistry department in
the
USA and to study something which would make me possible to collaborate with other members in the
institute. From 1962 I spent one year in the chemistry department of Indiana University (USA).
The main theme of the group I belonged to there was computational chemistry. I made a
poor computation person, and I regret now that I should have studied more seriously something about computation! This strange stay ended in a survey paper \cite{T5}.

Next year I went back to mathematics and spent one year in the Department
of Mathematics of California Institute of Technology with Professor W.A.J.~Luxemburg.
Our common interest was Riesz spaces (vector lattices). I was always nervous whether
I could write a paper during my stay. Happily near the end of my stay
I could finish a paper \cite{T9} on contractive projections in $L^p$ spaces.

The next long stay (1967--68) was in the Department of Mathematics of T\"ubingen University
(Germany) with Professor H.H.~Schaefer. Our common interest was Banach lattices. This was a pleasant period for me. I could write a paper \cite{T10}
on characterization of $L^p$ spaces in terms of existence of a positive projection
to every lattice subspace.

The third long stay (1975--76) was in University of Szeged (Hungary) with Professor
B.~Sz.-Nagy. Our common interest was Hilbert space operators. During my stay
I had a chance of writing a joint paper \cite{T13} with Dr. Z.~Ceau\c{s}escu, daughter of the late president of
Romania, and Prof. C.~Foia\c{s} on dilations of contractions.

In 1981 I spent two months in the Leningrad branch of the Institute of Mathematics of the Academy of Science, USSR by Professor N.~Nikolski. Our common interest was Hardy spaces.

In 1987 I spent two weeks in the Indian Statistical Institute, Delhi. Professor R.~Bhatia was my host. Our common interest was matrix theory.

Let me add some recent special visits. In 2011 summer, I visited Ferdowsi University of Mashhad, where I met an active group of young mathematicians conducted by Professor Moslehian, the second-named interviewer. On the same year I was
invited to the 60-th anniversary of Professor G.~Corach in Argentina, and in 2012 to the
60-th anniversary of Professor R.~Bhatia. I felt myself old enough as even those eminent
mathematicians had reached 60. It was nice for me to be able to invite some of those mathematicians mentioned above to Japan. Further Professor Moslehian joined in a conference in Japan.\\

\centerline{\includegraphics[keepaspectratio=true,scale=0.1]{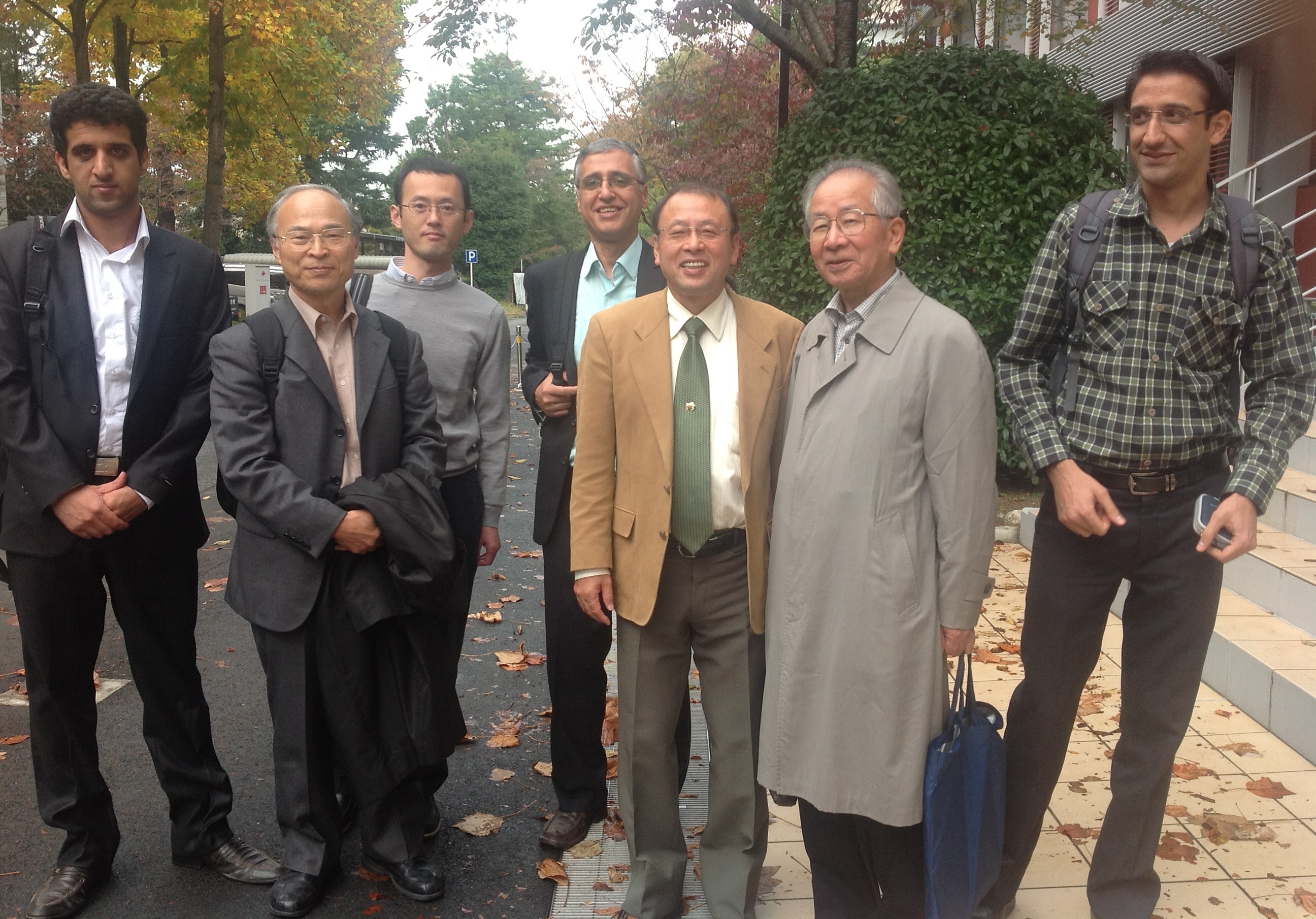}}
\vspace{.05in} \centerline{(H. Najafi (left), M. Uchiyama, M. Ito, M.S. Moslehian (back)}\vspace{.05in} \centerline{\newline M. Fujii, T. Ando (front), M. Kian (right), 2012)}
\vspace{0.2cm}

Professor Brualdi, the first-named reviewer, visited my Institute. As the editor-in-chief of the journal Linear Algebra and its Applications, Professor Brualdi is doing a great service to the researchers in this field. I am honored to have a
joint paper \cite{T25} with him.\\

\centerline{\includegraphics[keepaspectratio=true,scale=0.08]{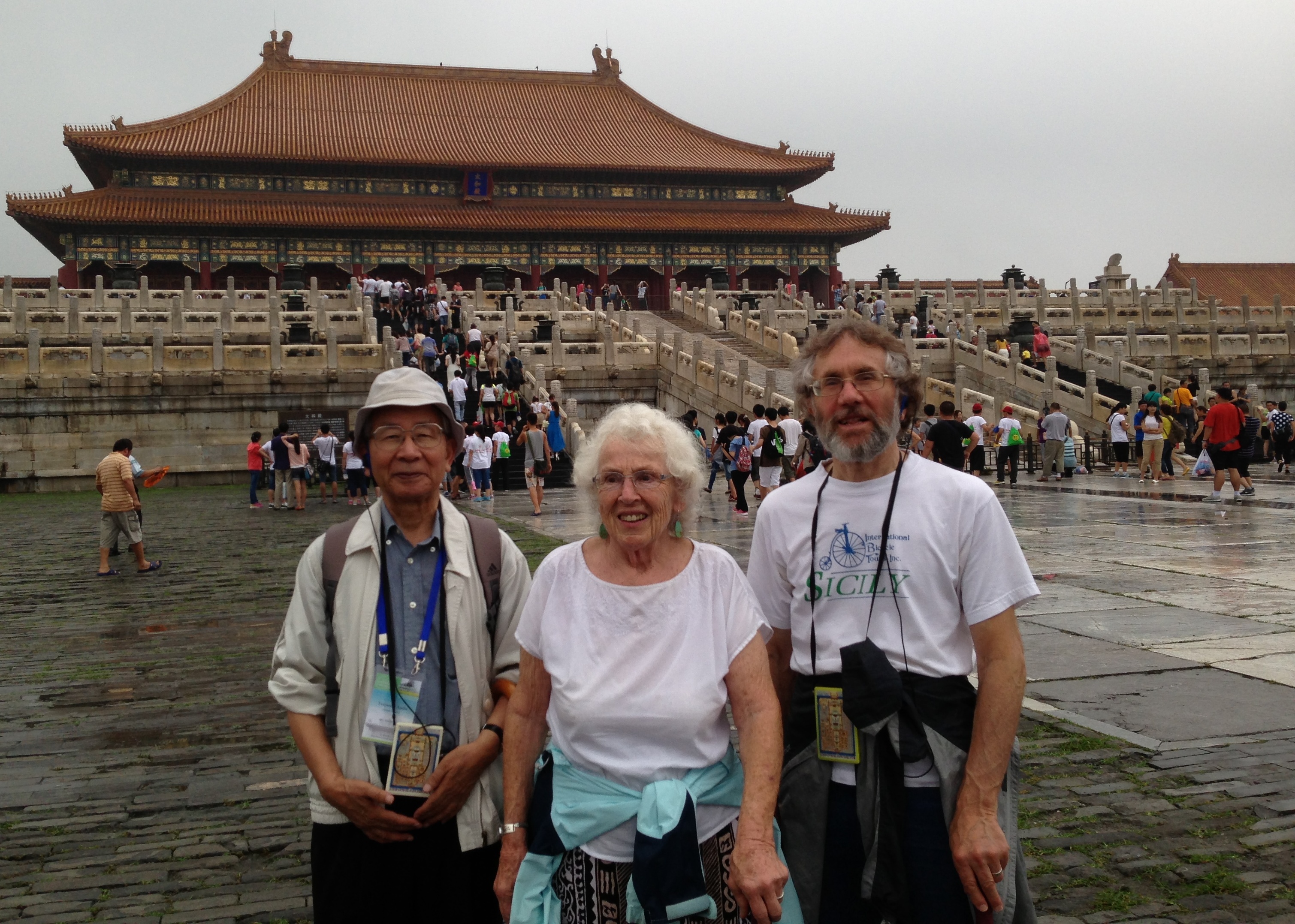}}
\vspace{.05in} \centerline{(T. Ando (left), M. Brualdi (middle), R.A. Brualdi (right), 2013)}

\vspace{.2in} $\bullet$ \textbf{Tell us how you do mathematics. Do you find writing mathematical
papers difficult?}\\

$\star$ There is no definite manner for me to do mathematics.
Writing a mathematical paper is a severe task. I feel uneasy whether this paper deserves publication, and I am always afraid that this paper would become my last paper.

\vspace{.2in} $\bullet$ \textbf{Who influenced you the most as a mathematician?}\\

$\star$
 I respected the achievements of my supervisor Professor Hidegoro Nakano though my personal contact with him or his guidance over me was small.
It is surprising to see that he established most of his fundamental contributions during the hard time of the war.
I inherited from him my life-long interest to the order structures of the objects.
The detailed account of his life and works can be found in a nice article by Prof. L.~Maligranda:
``Hidegoro Nakano (1909--1974)\ - \ on the centenary of his birth"
{\it Proceedings of the Third International Symposium on Banach and
Function Spaces III},
(Yokohama Publishers, 2011).
The other mathematician around me who gave influence over me is Dr.~Ichiro Amemiya, who came to
Hokkaido University with Prof. Nakano as a young lecturer.

Dr.~Amemiya was quite silent in contrast to talkative Nakano. Amemiya's ability in analyzing the essence of a problem and reformulating the problem was quite strong. I learned much from him about how to formulate problems. I wrote several papers with him. Let me cite two of them. One \cite{T7} is on the weak convergence of a random product of a finite number of selfadjoint contractions.
We proved that, given a finite number of selfadjoint contractions $T_1, T_2, \ldots, T_N$ in Hilbert space, every random product $T_{\sigma_1}T_{\sigma_2}\cdots T_{\sigma_n},$ where $\sigma_k \in
\{1,2, \ldots, N\}$ is a random choice, converges weakly to a limit, provided every $T_k$ appears infinitely often. At first I could not imagine the possibility of random choices.

\noindent Another paper \cite{T8} is on almost everywhere convergence
of the best $L^p$-approximants of a function from the subspaces of functions measurable with respect to a growing family of sub-$\sigma$ fields. This is a non-linear Martingale type theorem and we could also prove
the so-called maximal inequalities.

I admired the achievements of Mark Grigorievich Krein, a Soviet mathematician. This Jewish genius suffered very long segregation, yet established huge contributions in various branches of mathematics.

 I already mentioned the Krein--Rutman paper related to my thesis. I studied the following paper in detail:
M.G.~Krein, Theory of selfadjoint extensions of semi-bounded Hermitian operators
 and its applications I-II, {\it Mat. Sb.} 20(62) (1947) 431-495 and 21(63)(1947), 365-404.

Motivated by this first paper I wrote a paper \cite{T11} with K.~Nishio on extremal extensions
of positive symmetric operator. Also this paper of Krein
brought me to recognition of the order theoretic aspect of the Schur complement of a
positive semi-definite matrix (see \cite{T14} and \cite{T17}).

The monograph by F.G.~Gantmacher and M.G.~Krein,
``Oscillation matrices and small oscillations of mechanical systems",
seemed to appear even in 1941. But I could read its English translation only much later.
Here I knew the notion of total positivity of a matrix for the first time. I learned much about the anti-symmetric power of a matrix and wrote a lecture note on totally positive matrices (see \cite{T20}) on the basis of this monograph.

I studied several papers of M.G.~Krein (some with Ju. L.~Smuljan) on indefinite inner product spaces.
Those papers became the basis of my lecture note on the theory of Krein space. But my own contributions to the theory of
indefinite inner product spaces are not many but only \cite{T29} and \cite{T32}.

During my stay in Hungary, I studied in detail a series of works of M.G.~Krein with his
collaborators V.M.~Adamjan and D.Z.~Arov. Those papers became the basis of my lecture note on Hankel operators, and now are basic theories in Linear Systems Theory.

Thus the influence of M.G.~Krein to my works has been dominant. But I regret I did not study
 his more original and profound ideas based on classical analysis and his attitude
and ability toward concrete problems.

\vspace{.2in} $\bullet$ \textbf{What of your mathematical work do you like best? Tell us something
bout your mathematical work in general terms.}\\

$\star$ I started under the influence of Prof.~Nakano. In all of my works the standpoint of order-structure is dominant; in more detail,
order relations in vector lattices, matrix inequalities, norm inequalities, eigenvalue-majorization, operator-monotone functions etc..

I would like to summarize here the stream of my research. Let me start with some of my contributions in the field of Banach
lattices and related area. In this field nearest to the direction of Prof.~Nakano,
I tried to move in a different direction from him.

Most familiar examples of Banach lattices are $L^p$ spaces ($1 \leq p \leq \infty$). $L^2$ is specially a (real) Hilbert space in which every closed
subspace is the range of a projection of norm one. But the situation is
different {in other $L^p$s. There a typical example of a projection of norm one is a conditional expectation. In \cite{T9} I could characterize a projection of norm one in $L^p\ (p \not= 2)$ as a modified conditional expectation. The case
$p = 1$ was initiated by R.~Douglas.

In the converse direction I showed in \cite{T10} that if every closed lattice subspace of a Banach lattice becomes the range of a positive
projection, then the space itself should be an abstract $L^p$ space for some $p$.

The Orlicz space $L^{\Phi}$ with defining convex function $\Phi$ is also an important example of Banach lattices. The space $L^{\infty}$
can be seen as one of ``singular" Orlicz space with $\Phi(t) = 0$ or $\infty$ according as $0 \leq t \leq 1$ or $1 < t < \infty$.
Its dual is of $L^1$-type in the sense that the norm is additive on the cone of
non-negative elements. For a general Orlicz space its dual space admits
a unique order-orthogonal decomposition as the sum of order-continuous
functionals and that of order-singular ones. Nakano's interest seemed only in the order-continuous part. In \cite{T2} I showed that the (dual) norm is always of $L^1$-type on the order-singular part.

In connection with this order-continuous and singular decomposition
I found one interesting fact. On an abstract set provided
with a Boolean algebra of subsets one can consider the space of finitely additive measures. When the Boolean algebra is countably complete one can consider the (sub) space of countably additive measures. As in the case of $L^{\infty}$, each finitely additive measure is uniquely written
as a sum of a countably additive measure and a purely finitely additive one. The map ${\bf P}$ to the purely finitely additive
part is of course a projection of norm one. I could prove in \cite{T3} that
the projection ${\bf P}$ is also sequentially weak-continuous in the sense that
if a sequence of measures $\mu_n$ converge weakly to $\mu$ as $n \rightarrow \infty$ then ${\bf P}\mu_n $ converge weakly to ${\bf P}\mu$.
Afterward I recognized that this fact could follow from a theorem of A. Grothendieck.

A non-commutative (complex) version of $L^{\infty}$ is the von Neumann algebra,
which admits a unique predual. One of the analytic versions of $L^{\infty}$ is the Hardy space $H^{\infty}$. I showed in \cite{T15} that the Hardy space $H^1$ is the unique predual of $H^{\infty}$.

In the general theory of ordered Banach spaces it is usually easy
to determine, given some property for the space, what property the dual space should have. But the converse direction seemed unusually difficult. I challenged myself to a converse problem and showed in \cite{T4} that the dual
space has lattice-property if and only if the original space has the so-called
``interpolation-property", that is, if $a, b \leq c, d$ there is
an element $e$ such that $a, b \leq e \leq c, d$.

Next let me mention here something about my small contributions to the Hilbert spaces operators. I mentioned already about a theorem on convergence of random products of selfadjoint contractions.

Some day, in the seminar conducted by Prof.~Nakano, I was informed about a problem of
commuting unitary dilations. It asks, given a pair of commuting contractions $A, B$, whether there exists a commuting pair of unitary operators $U,V$ on a superspace such that
\[
\langle A^mB^nx|y\rangle = \langle U^mV^nx|y\rangle\quad {\rm ~for~all~}\ m, n,\ {\rm ~for~all~}\ x, y.
\]
I became aware of a small trick based on commutativity of $A, B$, which could give an affirmative answer \cite{T6}. I should develop this direction further, but
I could not.

Motivated by a result of C.A.~Berger, I had interest in getting a canonical form of a
numerical contraction, that is, \ $
w(T) := {\rm sup}\{|\langle Tx|x\rangle| :\ \|x\| = 1\} \leq 1.$

With a rather complicated recursive approach I arrived at a goal in \cite{T12}. Now it is known that there are several simpler ways to get such a canonical form. Once again afterward
 I returned back to a problem related to this norm $w(\cdot)$ and, together with K.~Okubo, using a technique of completely positive maps, we could characterize in \cite{T23} when, given a matrix $A$, the map $X \longmapsto X\circ A$ (Schur product) is contractive with respect to the norm $w(\cdot)$.

Most common order-relation for a pair $A, B$ of selfadjoint operators (or matrices) is the so-called L\"owner order, based on positive semidefiniteness $\geq 0$, that is,
$A \geq B$ means $A - B \geq 0.$ In contrast to the vector-lattice case, i.e. commutative case, a bounded family of selfadjoint operators does not admit its supremum (or infimum) in general. In the paper of M.G.~Krein mentioned above, I learned that for a block-matrix ${\bf A} = \begin{bmatrix}A_{11} & A_{12} \\ A_{21} & A_{22}\end{bmatrix} \geq 0$ the set
\[
\Big\{X \geq 0 ;
\begin{bmatrix}A_{11} & A_{12} \\ A_{21} & A_{22}\end{bmatrix}
 \geq \begin{bmatrix}X& 0 \\ 0 & 0\end{bmatrix}\Big\}
\]
admits the maximum element $A_{11} - A_{12}A_{22}^{-1}A_{21}$, which is just the Schur complement of ${\bf A}$ with respect to $A_{22}.$ Around the same period I encountered
 several results related to the parallel sum $A:B \overset{def}{=} (A^{-1} + B^{-1})^{-1}$ for $A, B > 0$ developed by a group of mathematicians under the influence of R.~Duffin as a matrix-version of parallel
connection of resistive networks in electrical network theory. Twice $A:B$ can be understood as the harmonic mean of $A$ and $B$. Then there came a question what is the geometric mean of $A$ and $B$.

In the space of $n\times n$ matrices the notion of positive linear maps is quite common. When restricted to the space of selfadjoint matrices there is also a natural non-linear map
induced via functional calculus by a real value function $f(t)$ defined on an interval $J \subset {\mathbb R}.$
When all eigenvalues of $A$ are contained in the interval, the functional calculus
 $f(A)$ is defined via the spectral decomposition of $A$. When $f(t) = \sum_{j=0}^k\alpha_jt^j$
is a polynomial, $f(A)$ is nothing but $\sum_{j=0}^k\alpha_jA^j. $ A natural question is
when the non-linear map $A \longmapsto f(A)$ preserves the order-relation. When it is
the case for all order $n$, the function $f(t)$ is said to be operator-monotone (increasing).\\

\centerline{\includegraphics[keepaspectratio=true,scale=0.17]{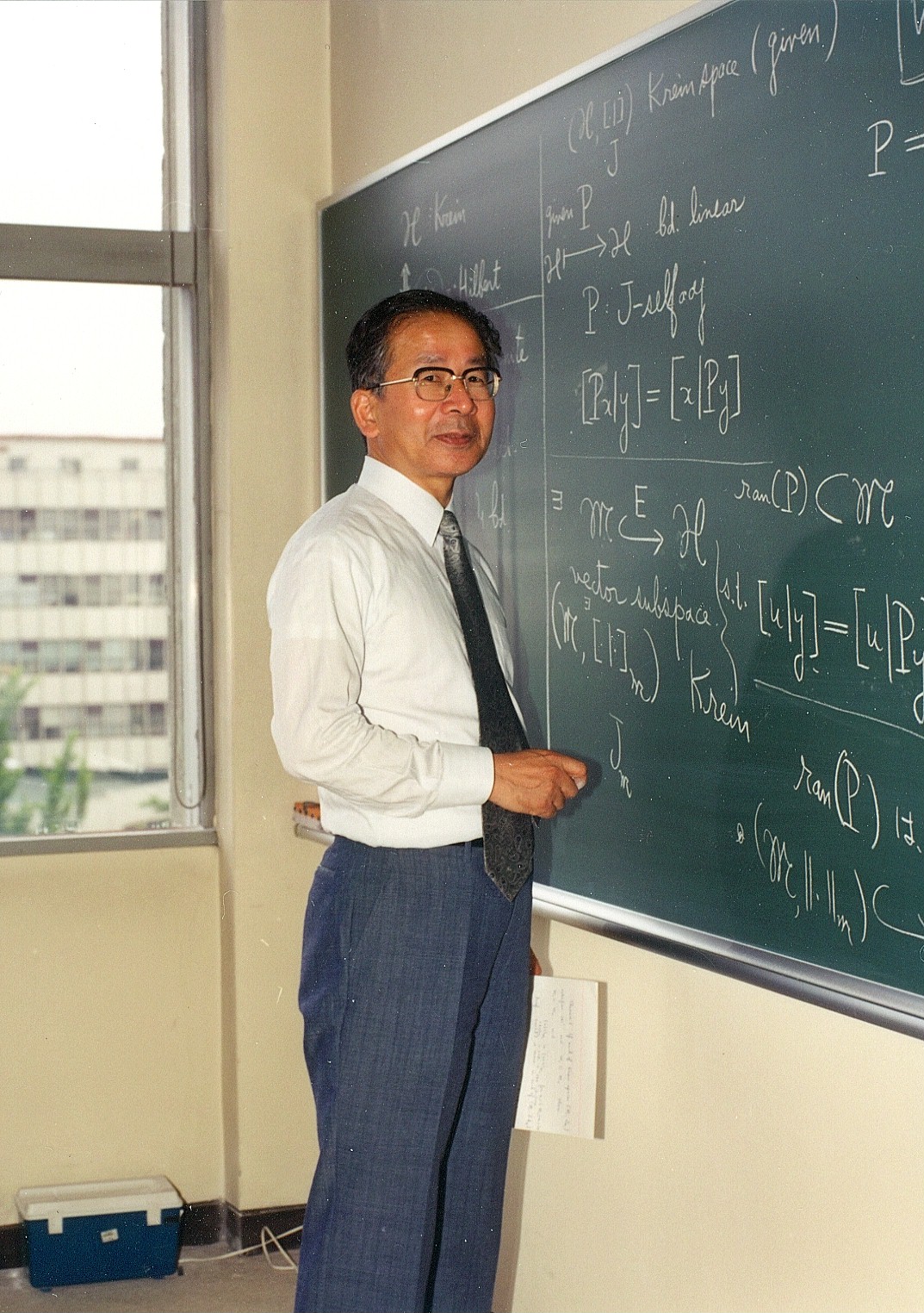}}
\vspace{.05in} \centerline{(T. Ando at a seminar in the Institute, around 1985)}
\vspace{0.5cm}

Surprisingly, already in 1933 K.~L\"owner succeeded in characterizing operator-monotonicity in term of positive semi-definiteness of the so-called L\"owner matrices.
Later a deep connection was discovered between operator-monotonicity and the Nevanlinna-Pick theory; the theory of function
which admits analytic continuation to the upper-half plane and maps the half-plane
into itself. Finally it turns out that when $f(t) > 0$ on $(0,\infty)$ is operator-monotone, it admits a (unique) integral representation
\[
f(t) = \alpha + \beta t + \int_{(0,\infty)}\frac{\lambda t}{t + \lambda}d\sigma(\lambda)
\]
where $\alpha, \beta \geq 0$ and $\sigma(\cdot)$ is a positive measure on $(0,\infty).$

F.~Kubo and I in \cite{T18} came to the conclusion
that
\[(A,B) \longmapsto A^{\frac{1}{2}}\cdot
f\Big(A^{-\frac{1}{2}}BA^{-\frac{1}{2}}\Big)\cdot A^{\frac{1}{2}}\]
 may give a nice definition
of a generalized mean and that every such mean can be represented as a generalized
 average of
 weighted arithmetic means $\alpha A + \beta B$ and weighted
harmonic means $(\alpha A):(\beta B)$ for $\alpha, \beta \geq 0$. Various properties of a mean can be derived from the properties
of weighted arithmetic and harmonic means. The most interesting is the geometric mean, written as $A\sharp B$, corresponding to the square-root function $f(t) = \sqrt{t}.$ The mean corresponding to the operator-monotone function $f(t) = t^{\alpha} \ (0 < \alpha < 1)$ is usually denoted by $A\sharp_{\alpha}B$. It is interesting to see that $A\sharp B$ can be characterized as
\[
A\sharp B = {\rm max}\Big\{X \geq 0 ; \begin{bmatrix}A & X \\ X & B\end{bmatrix} \geq 0\Big\}.
\]
The integral representations for operator monotone functions turned out quite useful
to produce new matrix inequalities. In \cite{T16} this method was applied to the case of tensor products of positive semi-definite matrices.

As mentioned before, the Schur product $A\circ B$ of $n\times n$ matrices $A, B$
is a principal submatrix of their tensor product $A\otimes B$, whose eigenvalues are simply all products of eigenvalues of $A$ and those of $B$. Therefore it is natural to expect that the eigenvalues or singular values of $A\circ B$ and those of $A$ and $B$,or even the eigenvalues or singular values of the usual product $AB$ are related in some manner.

It was realized in \cite{T19} with R.~Horn and Ch.R.~Johnson as an improvement of a known majorization relation for singular values
\[
\sum_{j=1}^k\sigma_j(A\circ B) \leq \sum_{j=1}^k\sigma_j(A)\cdot\sigma_j(B)\quad k = 1, \ldots, n
\]
and in \cite{T26} as a logarithmic majorization relation for eigenvalues
\[
 \prod_{j=k}^n\lambda_j(A\circ B) \geq \prod_{j=k}^n\lambda_j(AB) \quad k = 1,\ldots, n,
\quad {\rm for} \ A, B > 0
\]
where the singular values $\sigma_j(A)$ and eigenvalues $\lambda_j(A)$ of $A$, say, are indexed in non-increasing order respectively. (Notice that
though the usual product $AB$ is not selfadjoint all of its eigenvalues are real.)

In the space of $n\times n$ matrices most of useful norms $\|\cdot\|$ are unitarily invariant in the sense
\[
\|UAV\| = \|A\| \quad {\rm ~for~all~}\ A \quad {\rm ~for~all~} \ {\rm unitary}\ U, V.
\]
Every unitarily invariant norm $\|A\|$ is determined only by the singular values of $A$.
The Ky Fan dominance theorem says that majorization relation for singular values
\[
\sum_{j=1}^ks_j(B) \leq \sum_{j=1}^ks_j(A)\quad k = 1, \ldots, n
\]
implies the norm inequality $\|B\| \leq \|A\|$ for all unitarily invariant norms $\|\cdot\|.$

On this basis, I showed in \cite{T21} and with X.Zhan in \cite{T26} that for all positive operator-monotone function $f(t)$ on
$[0,\infty)$
\[
\|f(A) - f(B)\| \leq \|f(|A - B|)\| \,\, {\rm and}\,\,
\|f(A + B)\| \leq \|f(A) + f(B)\|\,\, {\rm ~for~all~} \ A, B \geq 0.
\]

Motivated by a result of H.~Araki, I together with F.~Hiai established in \cite{T24} a logarithmic-majorization relation for eigenvalues as, for $A, B > 0$
\[
\prod_{j=1}^k\lambda_j(A^r\sharp_{\alpha}B^r) \leq \prod_{j=1}^k\lambda_j (A\sharp_{\alpha}B)^r \ (k = 1,\ldots, n)\ {\rm ~for~all~} \
0 \leq \alpha \leq 1,\ {\rm and}\ r \geq 1.
\]
The essential part of this approach is in an implication relation among matrix inequalities,
which turns out to cover an interesting matrix inequality due to T.~Furuta.

\vspace{.2in} $\bullet$ \textbf{ What was your priority as a faculty member: Teaching or Research?}\\

$\star$
As I worked not in the mathematics department but in a research institute for electronics, my priority as a faculty member was surely research. Of course, as I was
associated to graduate education of the mathematics department, I gave a lecture
in most semesters.

Can you imagine the situations where you are the only mathematician in the whole institute and you should do only research ?

Some mathematics professor may say it would be wonderful if he is free from teaching duty. I don't agree with him. Suppose that you can devote full time to your research and
people outside expect remarkable results from you within certain years. There will be quite strong pressure upon you and you might find it difficult to bear.

Some teaching duty is good or necessary for your mental health. Every mathematician
should teach. It would be nice to be surrounded by talented students.

\vspace{.2in} $\bullet$ \textbf{When you are not doing mathematics, what do you do? What interests do
you have in addition to mathematics?}\\

$\star$
I am unhappy without a definite hobby. After my second retirement at the age of 70, I began to learn the Korean language to prevent me from senility. Between the Japanese language and the Korean language the grammars are almost similar, but daily life vocabularies are quite different. Advance in my learning is very slow, and situations in the future seem
hopeless!\\

\centerline{\includegraphics[keepaspectratio=true,scale=0.15]{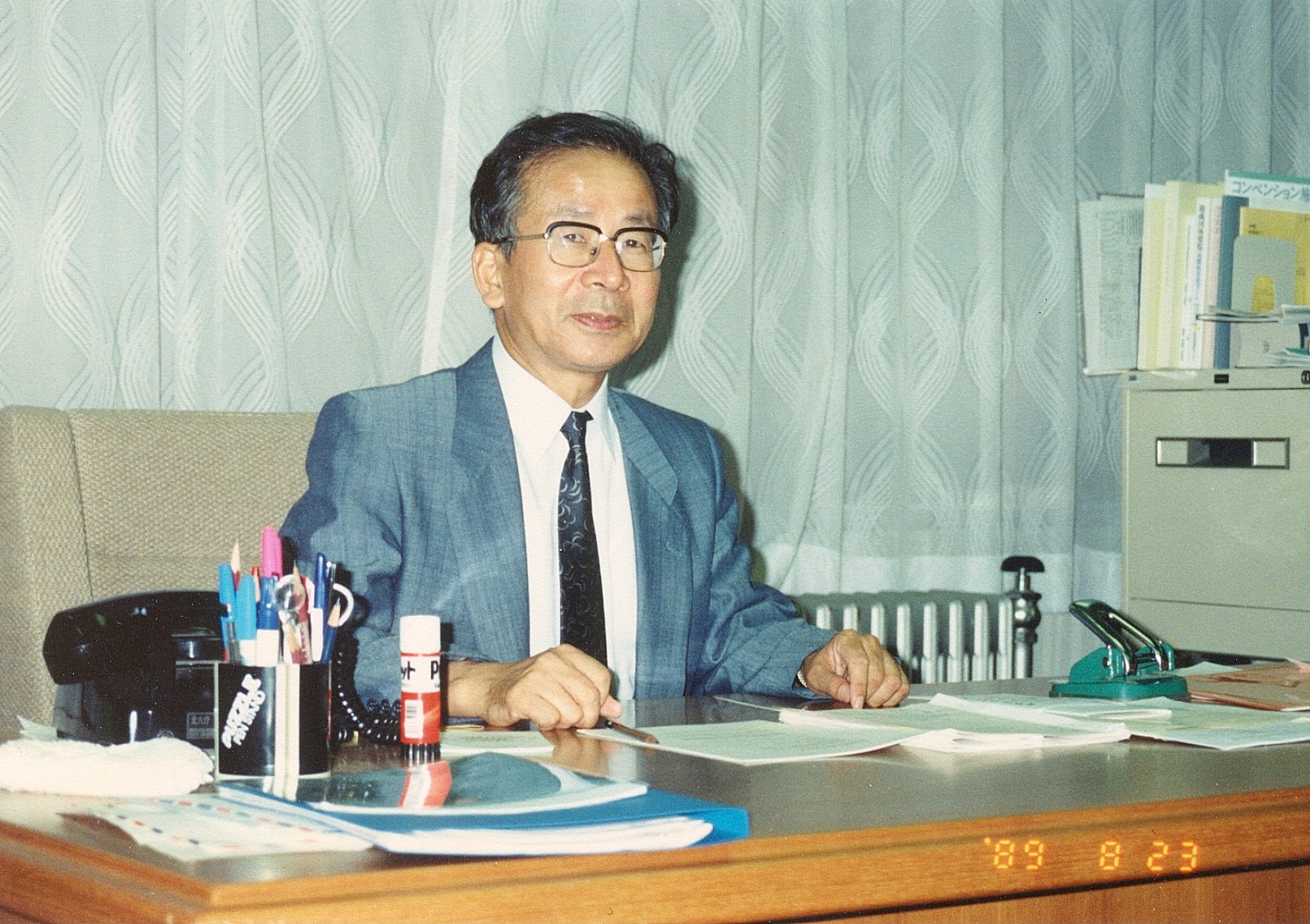}}
\vspace{.05in} \centerline{(T. Ando, at the director's room of the Institute, 1989)}

\vspace{.2in} $\bullet$ \textbf{How do you arrive at new ideas? Have examples/experimentation
been important in your work?}\\

$\star$ New ideas did not come suddenly, but only after many trial and errors.
To construct (counter) examples or to check examples was hard to me. I was not
skilled in computation.

\vspace{.2in} $\bullet$ \textbf{How do you think we can make mathematics more attractive for young
students? Is mathematics popular in Japan now?}\\

$\star$
As I was not much associated with undergraduate education of mathematics, let me say something about mathematics in high school. To encourage talented pupils, Japan has joined the International Mathematical Olympiad since 1990. The best record was the second in 2009. As plane geometry and combinatorics were not in the high school curricula, Japanese pupils seemed disadvantageous.

There is no nationwide contest of mathematics like the Et\"ov\"os contest in Hungary. Yet provincially, in Hokkaido, the capital of which is my city, there has been the Hokkaido High School Mathematics Contest for these 20 years.

\vspace{.2in} $\bullet$ \textbf{You wrote several books. Which one do you like more and why?}\\

$\star$
No, I have not written any monograph. Instead, I wrote several lecture notes and research reports:

\begin{itemize}
\item ``Topics on operator inequalities", Lecture Notes, Hokkaido Univ. 1979.
\item ``Linear operators on Krein spaces", Lecture Note, Hokkaido Univ. 1979.\\
\quad (The above two were accepted widely.)
\item ``Hankel operators I", Lecture Notes, Hokkaido Univ. 1983.
\item ``Reproducing kernel spaces and quadratic inequalities", Lecture Notes, Hokkaido Univ. 1987.
\item ``Totally positive matrices", Lecture Notes, Hokkaido Univ.1987. See \cite{T20}.
\item ``Majorization, doubly stochastic matrices and comparison of eigenvalues", 1987. See \cite{T22}.
\item ``Matrix quadratic equations", Lecture Notes, Hokkaido Univ. 1988.
\item ``De Branges spaces and analytic operator functions", Lecture Notes, Hokkaido Univ. 1990.
\item ``Completely positive matrices", Lecture Notes, Hokkaido Univ. 1991.
\item ``Operator-theoretic methods for matrix inequalities", Research Report, Hokusei Gakuen Univ. 1998.
\item Norms
and cones in tensor products of matrices", Research Report, Hokusei Gakuen Univ. 2001. See \cite{T28}.
\end{itemize}
I sent some of those lecture notes to Mathematical Review for review while some of them were published in the journal Linear Algebra and its Applications.

In addition, I wrote a survey article on Schur complements \cite{T31}.

\vspace{.2in} $\bullet$ \textbf{What is mathematical research, in your opinion? }\\

$\star$
It is difficult for me to answer. I can say nothing instructive about mathematics or true mathematical research. Each mathematician can do something of his ability. But about his attitude toward problems,
let me cite some words of Herman Weyl, a great mathematician, which I encountered in
 a book.

``But definite concrete problems were first conquered in their undivided complexity,
single handed by brute force, so to speak. Only afterwards the axiomaticians came along and stated: Instead of breaking in the door with all your might and bruising your hands, you should have constructed such and such a key of skill, and by it you would have been able to open the door quite smoothly. But they can construct the key only because they are able, after the breaking in was successful, to study the lock from within and without. Before you generalize, formalize, and axiomatize, there must be mathematical substance."

\vspace{.2in} $\bullet$ \textbf{ What is your opinion about the future of electronic publicatons?}\\

$\star$
Electronic publications will become a main trend. As to electronic journals, however, there may be a problem. Quick publications are welcome, but severe peer review should be kept.

Professor Moslehian is active enough to publish two electronic journals: Banach Journal of Mathematical Analysis and Annals of Functional Analysis. Both are getting growing reputation. I appreciate his efforts.

\vspace{.2in} $\bullet$ \textbf{Any final comments?}\\

$\star$
Perhaps some words of reflection:
In a small conference in honor of my retirement in 1995 I spoke about the order-structure of
tensor products of matrix spaces. (Some part of my talk is incorporated in the survey paper \cite{T28}.) But I could not go further.

At that time I did not recognize that rapid development in this direction had started among
physicists concerned with quantum information and computation. I am trying to
understand the stream of this development, hoping to make some contributions in this field.

Since the appearance of the paper of Ando--Li--Mathias \cite{T30} on geometric means for more than three positive definite matrices,
 many mathematicians tried to find what is the most suitable definition of such geometric mean. The direction based on some manifold theoretical observations seems most natural and productive. As the development to this direction is sometimes beyond my understanding, I am studying the background, hoping to make some contributions even in this direction.

\vspace{.2in} $\bullet$ \textbf{Thank you so much.}

\bigskip

\bibliographystyle{amsplain}

\begin{thebibliography}{50}


\bibitem{T7} I. Amemiya and T. Ando, \textit{Convergence of random products of contractions in Hilbert space}, Acta Sci. Math. (Szeged) \textbf{26} (1965), 239--244.

\bibitem{T1} T. Ando, \textit{Positive linear operators in semi-ordered linear spaces}, J. Fac. Sci. Hokkaido Univ. Ser. I \textbf{13} (1957), 214--228.

\bibitem{T2} T. Ando, \textit{Linear Functionals on Orlicz spaces}, Niew Arch. Wisk. (3) \textbf{8} (1960), 1--16.

\bibitem{T3} T. Ando, \textit{Convergent sequences of finitely additive measures}, Pacific J. Math. \textbf{11} (1961), 395--404.

\bibitem{T4} T. Ando, \textit{On fundamental properties of a Banach space with a cone}, Pacific J. Math. \textbf{12} (1962), 1163--1169.

\bibitem{T5} T. Ando, \textit{Properties of Fermion density matrices}, Rev. Modern Phys. \textbf{35} (1963), 690--702.

\bibitem{T6} T. Ando, \textit{On a pair of commutative contractions}, Acta Sci. Math. (Szeged) \textbf{24} (1963), 88--90.


\bibitem{T9} T. Ando, \textit{Contractive projections in $L_p$ spaces}, Pacific J. Math. \textbf{17} (1966), 391--405.

\bibitem{T10} T. Ando, \textit{Banachverb\"ande und positive Projektionen}, Math. Z. \textbf{109} (1969), 121--130.

\bibitem{T12} T. Ando, \textit{Structure of operators with numerical radius one}, Acta Sci. Math. (Szeged) \textbf{33} (1973), 11--15.


\bibitem{T14} T. Ando, \textit{Lebesgue-type decomposition of positive operators}, Acta Sci. Math. (Szeged) \textbf{38} (1976), 253--260.


\bibitem{T15} T. Ando, \textit{On the predual of $H^{\infty}$}, Comment. Math. \textbf{1} (1978), 33--40.

\bibitem{T16} T. Ando, \textit{Concavity of certain maps on positive definite matrices and applications to Hadamard products}, Linear Algebra Appl. \textbf{26} (1979), 173--186.

\bibitem{T17} T. Ando, \textit{Generalized Schur complements}, Linear Algebra Appl. \textbf{27} (1979), 173--186.


\bibitem{T20} T. Ando, \textit{Totally positive matrices}, Linear Algebra Appl. \textbf{90} (1987), 165--219.

\bibitem{T21} T. Ando, \textit{Comparison of norms $\|f(A) - f(B)\|$ and $\|f(|A - B|)\|$}, Math. Z. \textbf{197} (1988), 403--409.


\bibitem{T22} T. Ando, \textit{Majorization, doubly stochastic matrices and comparison of eigenvalues}, Linear Algebra Appl. \textbf{118} (1989), 163--248.

\bibitem{T26} T. Ando, \textit{Majorization relations for Hadamard products}, Linear Algebra Appl. \textbf{223/224} (1995), 283--295.

\bibitem{T28} T. Ando, \textit{Cones and norms in the tensor product of matrix spaces},
Linear Algebra Appl. \textbf{379} (2004), 3--41.

\bibitem{T29} T. Ando, \textit{L\"owner inequality of indefinite type}, Linear Algebra Appl. \textbf{385} (2004), 21--39.


\bibitem{T31} T. Ando, \textit{Schur complements and matrix inequalities, operator-theoretic approach}, pp.137--162 in ``The Schur complements and its applications" ed. F. Zhang, Springer, 2005.

\bibitem{T32} T. Ando, \textit{Projections in Krein spaces}, Linear Algebra Appl. \textbf{431} (2009), 2346--2355.


\bibitem{T8} T. Ando and I. Amemiya, \textit{Almost everywhere convergence of prediction sequences in $L_p (1 < p < \infty)$}, Z. Wahrscheinlichkeitstheorie und Verw. Gebiete \textbf{4} 1965 113--120 (1965).

\bibitem{T25} T. Ando and R. Brualdi, \textit{Sign--Central matrices,} Linear Algebra Appl. \textbf{208/209} (1994), 283--295.

\bibitem{T13} T. Ando, Z. Ceau\c{s}escu and C. Foia\c{s}, \textit{On intertwining dilations. II}, Acta Sci. Math. (Szeged) \textbf{38} (1976), 3--14.

\bibitem{T24} T. Ando and F. Hiai, \textit{Log majorization and complementary Golden-Thompson type inequalities}, Linear Algebra Appl. \textbf{194} (1994), 113--131.

\bibitem{T19} T. Ando, R. Horn and Ch.R. Johnson, \textit{The singular values of a Hadamard product: a basic inequality}, Linear Multilinear Algebra \textbf{21} (1987), 345--365.

\bibitem{T30} T. Ando, C.K. Li and R. Mathias, \textit{Geometric means}, Linear Algebra Appl. \textbf{385} (2004), 305--334.

\bibitem{T11} T. Ando and K. Nishio, \textit{Positive selfadjoint extensions of positive symmetric operators}, Tohoku Math. J. (2) \textbf{22} (1970), 65--75.

\bibitem{T23} T. Ando and K. Okubo, \textit{Induced norms of the Schur multiplier operator},
Linear Algebra Appl. \textbf{147} (1991), 181--199.

\bibitem{T27} T. Ando and X. Zhan, \textit{Norm inequalities related to operator monotone functions}, Math. Ann. \textbf{315} (1999), 771--780.

\bibitem{T18} F. Kubo and T. Ando, \textit{Means of positive linear operators}, Math. Ann. \textbf{246} (1979/80), 205--224.


\end{thebibliography}

\end{document}